\documentclass{amsart}
\usepackage{amssymb}
\begin{document}

% Redefine "plain" pagestyle
%%%\makeatletter	   % `@' is now a normal "letter' for LaTeX
%%%\renewcommand{\ps@plain}{%
     %%%\renewcommand{\@oddhead}{\textrm{Your Header}\hfil\textrm{\thepage}}% 
     %%%\renewcommand{\@evenhead}{\@oddhead}%
     %%%\renewcommand{\@oddfoot}{}% empty footer
     %%%\renewcommand{\@evenfoot}{\@oddfoot}}
%%%\makeatother     % `@' is restored as a "non-letter" character

\def\bp{\mathbb P}
\def\bpn{\mathbb P^n}
\def\bc{\mathbb C}
\def\ix{I(X)}
\def\iy{I(Y)}
\def\iz{I(Z)}
\def\so{S_1\otimes}
\def\sx{S(X)}
\def\bwdg{\ifmmode \hbox{$\bigwedge$} \else \bigwedge \fi}
	%% definition created by A. Vogel 2/25/98 
\def\bwdgl{\ifmmode \hbox{$\bigwedge\nolimits$} \else      
   \bigwedge\nolimits \fi }

%Declarations
\theoremstyle{plain}
\newtheorem{thm}{Theorem}[section]
\newtheorem{cor}[thm]{Corollary}
\newtheorem{lem}[thm]{Lemma}
\newtheorem{clm}[thm]{Claim}
\newtheorem{quest}[thm]{Question}
\theoremstyle{definition}
\newtheorem{defn}[thm]{Definition}

\newtheorem{rem}[thm]{Remark}

\newcommand{\kindalong}{\hbox{$\hbox to .35in{\rightarrowfill}$}  }
\newcommand{\reallylong}{\hbox{$\hbox to .5in{\rightarrowfill}$}  }

%Commands
%%%\errorcontextlines=0
\numberwithin{equation}{section}

\title{Resolutions of Subsets of Finite Sets of Points in Projective Space}         
\author{Steven P. Diaz}  
\address{(S. P. Diaz) Department of Mathematics, Syracuse University\\
                  Syracuse, N.Y. 13244  United States of America}
\email{spdiaz@mailbox.syr.edu}    
\author{Anthony V. Geramita}  
\address{(A. V. Geramita) Department of Mathematics and Statistics, Queen's 
University\\
        Kingston, Ontario, Canada K7L 3N6           }
\email{tony@mast.queensu.ca}   
\address{(A. V. Geramita) Dipartimento Di Matematica, Universit\'a Di Genova\\
          Genova, Italia}
\email{geramita@dima.unige.it}
\author{Juan C. Migliore}  
\address{(J. C. Migliore) Department of Mathematics, University of Notre Dame\\
          Notre Dame, Indiana  46556  United States of America }
\email{migliore.1@nd.edu} 
%%%\date{}  

\begin{abstract}
Given a finite set, $X$, of points in projective space for which the Hilbert
function is known, a standard result says that there exists a subset of this
finite set whose Hilbert function is ``as big as possible'' inside $X$.  Given
a finite set of points in projective space for which the minimal free resolution
of its homogeneous ideal is known, what can be said about possible resolutions
of ideals of subsets of this finite set?  We first
give a maximal rank type description of the most generic possible resolution of
a subset.  Then we show that this generic resolution is not always achieved, by
incorporating an example of Eisenbud and Popescu.  However, we show that it
{\em is} achieved for sets of points in projective two space: given any finite
set of points in projective two space for which the minimal free resolution is
known, there must exist a subset having the predicted resolution.
\end{abstract}
        
\maketitle

\pagestyle{plain}

%**************************************************** 

\section{Introduction}       
We work over an algebraically closed field, $k$.  Let  $X=\{P_1, \dots,
P_d\}$ be a finite set of $d$ distinct points in projective $n$-space over
$k$, $\bpn$.  Associated to $X$ we have its homogeneous ideal
$\ix\subset k[x_0, \dots, x_n] = S$ and its homogeneous coordinate ring $\sx
= S/\ix$.  A fundamental invariant of $X$ is its Hilbert function, $h_X$,
defined to be the Hilbert function of $\sx$:
\begin{equation*}
h_X(t) = \dim_{k}\sx_t.
\end{equation*}
Lacking some uniformity property such as the Uniform Position Property (UPP),
the subsets of $X$ of fixed cardinality may have different Hilbert functions. 
Given this, it is somewhat surprising at first glance that there is always at
least one subset with a predetermined Hilbert function.  Indeed, one of the
fundamental results about Hilbert functions of subsets of $X$ is  the following:

\begin{lem}\label{lem:1}
Fix an integer $e$, $1\le e < d$. Then there exists a subset $Y$ of $X$ of 
exactly $e$ points such that
\begin{equation*}
h_Y(t) = \min\{h_X(t), e\}.
\end{equation*}
\end{lem}
\begin{proof}
This is very well known. See for instance \cite{GMR}, Lemma~2.5~(c).  See also
Remark \ref{wrong-approach}.
\end{proof}

The Hilbert function is a very coarse measure of the properties of $X$.  
Related finer measures that are often studied are the graded Betti numbers
and twists of the minimal graded free resolution of $\sx$ or equivalently
$\ix$.
\begin{equation}\label{eq:1.2}
0\to F_{n-1}\to \dots \to F_1 \to F_0 \to \ix \to 0
\end{equation}
where $F_i = \oplus_{j=1}^{r_i}(S(-\gamma_{ij}))^{\alpha_{ij}}$. The 
$\gamma_{ij}$ are the twists and the $\alpha_{ij}$ are the graded Betti
numbers.  Since Lemma~\ref{lem:1} is so useful for Hilbert functions, one may
wonder whether there is a corresponding result for resolutions.

In section 2 we state a natural first guess at a possible generalization of 
Lemma~\ref{lem:1} to resolutions.  The guess is stated in terms of Koszul
cohomology.  The basic idea of the guess is that at least one subset of $X$
of each cardinality should behave as generically as possible subject to some
obvious constraints imposed by being a subset of $X$.  The guess is very similar
to the Minimal Resolution Conjecture of Lorenzini \cite{L2} except that the
Minimal Resolution Conjecture does not deal with subsets.

In section 3 we show that the guess of section 2 is incorrect.  In fact a 
counterexample to the Minimal Resolution Conjecture provided in \cite{EP} is
used to construct a counterexample to the guess.  While the Minimal
Resolution Conjecture did not turn out to be true in full generality, it is
true in many cases and is perhaps a good first start at understanding the
true situation.  (The end of the introduction to \cite{HS} contains a good
list of references to results about the Minimal Resolution Conjecture.)  One
might still hope that the guess of section 2 would behave similarly.  We make
this hope more precise with some questions at the end of section 3.

In sections 4 and 5 we answer these questions for $\bp^2$ by showing that
the guess is true for sets of points in $\bp^2$ (a place where the Minimal
Resolution Conjecture of Lorenzini is also known to be true).  A variety of
tools are used to carry this out.  We divide the problem into four cases,
depending on the number of minimal generators of
$I(X)$ in the maximum possible degree.  The three easiest of these cases are 
treated in section 4.  The most difficult is the case where $I(X)$ has two
minimal generators in this degree, and this case is treated in section 5.  Here
we combine liaison theory with a careful study of certain sections of a certain
twist of $\Omega_{\bp^2}^1$, the sheaf of differential one-forms on $\bp^2$.

%**************************************************** 

\section{A First Guess}
We first recall briefly how the graded Betti numbers of an ideal may be 
computed using Koszul cohomology; see \cite{G} section 1 for more details. 
One makes a complex
\begin{equation}\label{eq:2.1}
\begin{array}{rl}
\dots \to \bwdgl^{p+1}\so\ix_{q-1} & \stackrel{d_{p+1,q-1}}{\reallylong}
\bwdgl^p\so\ix_q \stackrel{d_{p, q}}{\kindalong} \bwdgl^{p-1}\so\ix_{q+1}
\\ 
& \stackrel{d_{p-1,q+1}}{\reallylong} \bwdgl^{p-2}\so\ix_{q+2} \to \dots 
\end{array}
\end{equation}
where $d_{p, q}(l_1\wedge l_2\wedge \dots \wedge l_p\otimes f) = 
\sum_{i=1}^{p}(-1)^{p-i}l_1\wedge \dots \wedge l_{i-1}\wedge l_{i+1}\wedge
\dots \wedge l_p\otimes l_if.$

In the resolution \eqref{eq:1.2} the exponent of $S(-(p+q))$ in $F_p$ is the 
dimension, as a vector space over $k$, of the cohomology group
\begin{equation*}
\frac{\ker d_{p, q}}{\text{im } d_{p+1, q-1}}.
\end{equation*}
Of course an exponent of $0$ means that $S(-(p+q))$ does not appear.  One 
certainly knows the dimensions of the vector spaces $\bigwedge^iS_1$.  If one
also knew the Hilbert function of $X$ and thus the dimensions of the vector
spaces $\ix_j$, then to compute all the graded Betti numbers and twists for
$\ix$ it would be sufficient to know the ranks of all the maps $d_{i, j}$. 
Thus, our guess will combine Lemma \ref{lem:1} with a guess about these ranks.

As before $X = \{P_1, \dots, P_d\}$ consists of $d$ distinct points.  We 
assume that the Hilbert function and resolution of $X$ are known.  We fix an
integer $1\le e < d$.  We guess that there should exist a subset $Y$ of $X$
of exactly $e$ points such that the graded Betti numbers and twists of the
graded minimal free resolution of $\iy$ are determined as follows in (a) and
(b).
\begin{itemize}
\item[(a)] The Hilbert function of $Y$ is as in Lemma~\ref{lem:1}.  Since 
$Y\subset X$, for each $i$, $\ix_i\subset \iy_i$ and we may compare the
complex (\ref{eq:2.1}) for $\ix$ and the corresponding one for $\iy$ by the
following commutative diagram.
\begin{alignat*}{2}
\dots \to\, \bwdgl^p&\so \ix_q\,& \stackrel{d_{p, q}}{\kindalong} 
\bwdgl^{p-1}&\so\ix_{q+1}\, \to\dots \\ &\quad\bigcap  &&\quad\bigcap  \\
\dots \to\, \bwdgl^p&\so \iy_q\,& \stackrel{e_{p, q}}{\kindalong}
\bwdgl^{p-1}&\so \iy_{q+1}\, \to\dots
\end{alignat*}
Assuming (a), we know the dimensions of all the $\iy_j$.  We then guess that 
the ranks of the $e_{i, j}$ will be as follows
\item[(b)] Of course $e_{0, p+q}$ is the zero map.  Having determined the 
rank of $e_{i, p+q-i}$, the rank of $e_{i+1, p+q-i-1}$ is as large as
possible subject to the two constraints:
	\begin{itemize}
	\item[(i)] $\ker e_{i+1, p+q-i-1}$ must contain $\ker d_{i+1, p+q-i-1}$
	\item[(ii)] $\text{im } e_{i+1, p+q-i-1}$ must be contained in $\ker e_{i, 
p+q-i}$.
	\end{itemize}
In other words rank $e_{i+1, p+q-i-1}$ is the smaller of
	\begin{itemize}
	\item[(i$'$)] $\dim\bigwedge^{i+1}\so\iy_{p+q-i-1} - \dim \ker d_{i+1, 
p+q-i-1}$  and
	\item[(ii$'$)] $\dim \ker e_{i, p+q-i}$.
	\end{itemize}
\end{itemize}

%**************************************************** 

\section{A Counter-Example to the First Guess}
As mentioned in the introduction this guess is similar to the Minimal 
Resolution Conjecture of Lorenzini.  They are both Maximal Rank Conjectures
in that they conjecture that certain vector space maps have ranks as large as
possible.  It is not surprising, therefore, that one can construct a
counterexample to the guess out of a counterexample to the Minimal Resolution
Conjecture.

In the introduction to \cite{EP} they point out that for 11 general points 
in $\bp^6$ the Minimal Resolution Conjecture predicts the resolution to be
\begin{align*}
0\to S(-8)^4 \to S(-7)^{18} \to S(-6)^{25}\oplus &S(-5)^{4} \to S(-4)^{45} 
\to \\ &S(-3)^{46} \to S(-2)^{17} \to I \to 0
\end{align*}
whereas the actual resolution is
\begin{align*}
0\to S(-8)^4 \to S(-7)^{18} \to S(-6)^{25} \oplus & \ S(-5)^5 \to S(-5)^1 
\oplus S(-4)^{45} \to \\ 
&S(-3)^{46} \to S(-2)^{17} \to I \to 0.
\end{align*}

 From \cite{L1} section 3 or \cite{L2} section 3 we know that the Minimal 
Resolution Conjecture is true for 22 general points in $\bp^6$.  Thus one can
work out that the resolution of 22 general points in $\bp^6$ is
\begin{align*}
0\to S(-8)^{15} \to S(-7)^{84} \to S(-6)^{190} \to &S(-5)^{216}  \to
S(-4)^{120} \to \\ &S(-2)^{6} \oplus S(-3)^{20} \to I \to 0.
\end{align*}
Any subset of 11 points of 22 general points is a set of 11 general points.  
Let us see what the guess predicts as the resolution of 11 points contained
in 22 general points.  The crucial term is $S(-5)$ so we compute only that. 
The relevant Koszul complex to look at is
\begin{equation}\label{eq:3.1}
0 \to\bwdgl^3\so I_2\,\stackrel{d_{3, 2}}{\kindalong} \bwdgl^2\so 
I_3\,\stackrel{d_{2, 3}}{\kindalong} \bwdgl^1\so
I_4\,\stackrel{d_{1, 4}}{\kindalong} \bwdgl^0\so I_5 \to 0.
\end{equation}

When $I$ is the ideal of 22 general points, using that these points  have 
generic Hilbert function $1, 7, 22, 22, \dots$, one easily computes the
dimensions in \eqref{eq:3.1} as
\begin{equation*}
0 \to 210 \, \stackrel{d_{3, 2}}{\kindalong} \, 1302 \, 
\stackrel{d_{2, 3}}{\kindalong}\, 1316\, \stackrel{d_{1,4}}{\kindalong} \, 440
\to 0.
\end{equation*}
To get the known resolution we must then have
\begin{alignat*}{2}
\text{rank }d_{1, 4} &= 440 &\qquad \dim\ker d_{1, 4} &= 1316 - 440 = 876\\
\text{rank }d_{2, 3} &= 876 &\qquad \dim\ker d_{2, 3} &= 1302 - 876 = 426\\
\text{rank }d_{3, 2} &= 210 &\qquad \dim\ker d_{3, 2} &= \;\; 210 - 210 = 0
\end{alignat*}
so that $\dim \frac{\ker d_{2, 3}}{\text{im } d_{3, 2}} = 426 - 210 = 216$.

When $I$ is the ideal of 11 general points, one again easily computes the 
dimensions in \eqref{eq:3.1} as
\begin{equation*}
0 \to 595 \, \stackrel{e_{3, 2}}{\kindalong}\, 1533 \,
\stackrel{e_{2, 3}}{\kindalong}\, 1393 \, \stackrel{e_{1,
4}}{\kindalong}\, 451 \to 0
\end{equation*}
Applying the guess we get
\begin{alignat*}{2}
\text{rank }e_{1, 4} &= 451 &\qquad \dim\ker e_{1, 4} &= 1393-451 = 942\\
\text{rank }e_{2, 3} &= 942 &\qquad \dim\ker e_{2, 3} &= 1533 - 942 = 591\\
\text{rank }e_{3, 2} &= 591 &\qquad \dim\ker e_{3, 2} &= \;\; 595 - 591 = 4
\end{alignat*}
so that $\dim \frac{\ker e_{2, 3}}{\text{im } e_{3, 2}} = 591 - 591 = 0$ and 
$\dim \frac{\ker e_{3, 2}}{\text{im } e_{4, 1}} = 4 - 0 = 4$.  That is, the
guess predicts the same resolution for 11 general points in $\bp^6$ as the
Minimal Resolution Conjecture, which is wrong.

The guess does give some restrictions on what resolutions of subsets of $X$ 
can be.  Conditions (i) and (ii) must always be satisfied, but the rank could
be smaller than this upper bound.  Also, because of the inductive nature of
the upper bounds, once one $e_{i, j}$ fails to achieve the upper bound, the
upper bounds on $e_{i+s, j-s}$ for $s\ge 1$ can change.

The above calculations give some preliminary evidence that the following
question may have an affirmative answer.  See also Remark \ref{MRC-in-P2}.

\begin{quest}
When $X$ is a general set of $d$ points in projective space, does the guess for
a subset of  $e < d$ points of $X$ always give the same graded Betti numbers as
those given by the Minimal Resolution Conjecture for a general set of $e$
points?
\end{quest}

Since the known counter-examples to the Minimal Resolution conjecture go wrong
in the ``middle'' of the resolution, it may well be that parts of
the Minimal Resolution Conjecture are always true.  In particular, it is known
that the Cohen-Macaulay Type Conjecture is true (\cite{trung-valla},
\cite{lauze}), and one naturally wonders if the Ideal Generation Conjecture is
true.  This leads to the second question:

\begin{quest} \label{q2}
Is the guess true at least at the ends of the resolution?  In particular,
given any finite set of points in $\bp^n$, is there always a subset with the
predicted minimal generators and the predicted Cohen-Macaulay type? 
\end{quest}

Note that the guess does not assume that we have a general set of points, or
even that we have some sort of uniformity!  The next two sections show that
Question \ref{q2} has an affirmative answer for subsets of $\bp^2$.

%%%%%%%%%%%%%%%%%%%%%%%%%%%%%%%%%%%%%%%%%%%%%%%%%%%%%%%%%%%%%%%%%%%%%%%%%

\section{The Subset Resolution Theorem for points in $\bp^2$}

We now restrict our attention to points in $\bp^2$.  Let $X = \{P_1, \dots, 
P_d\}$ be a set of $d$ distinct points, with homogeneous ideal $I =
\ix\subset k[X_0, X_1, X_2] = S$.  At first glance the Koszul complex would
seem to involve sequences of the form
\begin{align*}
0\to\bwdgl^3\so I_{s-2}\, \stackrel{d_{3, s-2}}{\kindalong} 
\bwdgl^2\so I_{s-1}\, &\stackrel{d_{2, s-1}}{\kindalong} \bwdgl^1\so
I_s\,\\ &\stackrel{d_{1, s}}{\kindalong} \bwdgl^0\so I_{s+1}\, \to 0
\end{align*}
However, we know that the graded minimal free resolution for $I$ has only 
two terms, so we must have that $d_{3, s-2}$ is injective and $\ker d_{2,
s-1} = \text{im } d_{3, s-2}$.  The only thing in question is the rank of the
map $d_{1, s}$.  This is just the multiplication map
\begin{align*}
\mu_s : \so & I_s \to I_{s+1}\\
L\otimes & F \mapsto LF.
\end{align*}

\begin{defn}If $Y \subset X$ has the Hilbert function given in Lemma
\ref{lem:1}, we will say that it has {\em truncated Hilbert function}.
\end{defn}

For any subset $Z\subset X$ and any positive integer $s$, we have a 
commutative diagram
\begin{equation}\label{eq:4.5}
\begin{array}{ccc}
\so\ix_s & \stackrel{\mu_{s, X}}{\kindalong} & \ix_{s+1} \\
 \bigcap & &\bigcap\\
\so\iz_s & \stackrel{\mu_{s, Z}}{\kindalong} &\iz_{s+1}
\end{array}
\end{equation}

The subset resolution guess for points in $\bp^2$ then becomes the 
following.
\begin{thm}\label{thm:4.6}
Let $X$ be a reduced set of $d$ points in $\bp^2$.  Fix an integer $m$, $1\le
m<d$.  Then there exists  a subset $Z \subset X$ of cardinality $m$ and with
truncated Hilbert function, as given in Lemma \ref{lem:1}, and such that
for all positive integers $s$
\begin{align*}
\text{rank } \mu_{s, Z} = \min\{&\dim \iz_{s+1}, \text{rank } \mu_{s, X} +\\
&\dim \so\iz_s - \dim \so \ix_s\}.
\end{align*}
\end{thm}

\begin{proof}
First observe that if $X$ is contained in a line then $X$ and all its subsets 
are complete intersections.  The resolution of a complete intersection is
well known.  We let the reader check the theorem in this case.  Thus we
may assume that $X$ is not contained in a line.

Now we show that it is enough to prove the theorem for $m=d-1$.  To do this, it
is enough to show the following.  Let $Z \subset X$ be a subset consisting of
$m=m_0$ points, such that $Z$ has truncated Hilbert function and $\mu_{s,Z}$
has the predicted rank, for any $s$.  Assume that there is a subset $Z_1
\subset Z$ consisting of $m_0-1$ points such that $Z_1$ has truncated Hilbert
function, and such that for all $s$, the rank of $\mu_{s,Z_1}$ is what is
predicted in the theorem if we take $X=Z$ and $d=m_0$.  Then we have to show
that this is the same rank that is predicted by the theorem if we had taken
$X=X$ and $Z = Z_1$.

The fact that $Z$ has the predicted rank says that for all $s$,  either $\mu_{s, Z}$
is surjective or $\ker \mu_{s, Z} = \ker \mu_{s, X}$.  By our assumption on
$Z_1$, we get that for all $s$ either $\mu_{s, Z_1}$ is surjective or $\ker
\mu_{s, Z_1} = \ker \mu_{s, Z}$.  For those $s$ with $\mu_{s, Z_1}$ surjective
we are done.  For those $s$ with $\ker\mu_{s, Z_1} = \ker\mu_{s, Z} =
\ker\mu_{s, X}$ we are done.  This leaves only those $s$ for which 
$\ker\mu_{s, Z_1} = \ker\mu_{s, Z}$ but $\ker\mu_{s, Z}\ne \ker\mu_{s, X}$. 
But if $\ker\mu_{s, Z}\ne \ker\mu_{s, X}$ then $\mu_{s, Z}$ is surjective. 
Furthermore, the Hilbert function of $X$ in degree $s$ is different from that
of $Z$ in degree $s$.  Since $Z$ has truncated Hilbert function, this says that
$Z$ imposes $m_0$ independent conditions on forms of degree $s$, and $Z_1$
imposes $m_0 -1$ independent conditions on forms of degree $s$.  Hence $\mu_{s,
Z_1}$ is surjective.  Thus it is enough to prove the theorem for the case
$m=d-1$.

Part of the proof will be by induction on $d$. The reader can easily get this
induction argument started by checking directly that the theorem is true for
small values of $d$.  

Let $l$ be the smallest positive integer such that $X$ imposes $d$ conditions 
on forms of degree $l$.  Then $\ix$ is generated in degrees less than or
equal to $l+1$, \cite{DGM} Prop.~3.7.  Let $Z\subset X$ be any subset with
$d-1$ points.  Then $\iz$ is also generated
in degrees less than or equal to $l+1$.  Thus for $s\ge l+1$ the multiplication
map $\mu_{s,Z}$ is surjective and thus satisfies the conclusion of the
proposition.  

Next consider $s\le l-1$.  {\em From now on, unless specified otherwise, we
assume that our subset $Z$ has truncated Hilbert function}.  Then $X$ imposes at
most $d-1$ conditions on forms of degree $s$. By Lemma~\ref{lem:1}, Z imposes
the smaller of $d-1$ and the number of conditions imposed by $X$ on
forms of degree $s$.  Thus, $\ix_s = \iz_s$.  This says that
$\mu_{s,X}$ and $\mu_{s,Z}$ are the same map.  Certainly the conclusion of
the proposition follows in this case.  We are only left to consider the case
$s=l$.  We have to show that {\em among subsets with cardinality $d-1$ and with
truncated Hilbert function}, we can find one with the right number of minimal
generators in degree $l+1$.

Consider the diagram \eqref{eq:4.5} with $s=l$.  Regardless of whether or not
$Z$ has truncated Hilbert function,  $\ix_l$ has codimension one  in $\iz_l$,
and similarly for $l+1$.  Let $F_1, \dots, F_t$ be a basis for
$\ix_l$ and let $G$ be a form of degree $l$ in $\iz_l - \ix_l$, so that $F_1,
\dots, F_t$, $G$ is a basis for $\iz_l$.  The proof breaks down into four
cases according to the codimension of the image of $\mu_{l,X}$ in
$\ix_{l+1}$, in other words, the number of generators $\ix$ needs in degree
$l+1$.  Note that if $\ix_l$ is zero dimensional then $\iz_l$ is one
dimensional, so $\mu_{l,Z}$ is injective.  We may assume that $\ix_l$ has
positive dimension.

Case 1.
$\ix$ is generated in degrees $\le l$.  This says that $\mu_{l,X}$ is 
surjective.  We wish to show that $\mu_{l,Z}$ is also surjective, for any $Z$
(hence in particular one with truncated Hilbert function).  Since
$\ix_{l+1}$ has codimension one in $\iz_{l+1}$ we simply need to find a
single form in the image of $\mu_{l,Z}$ not in the image of $\mu_{l,X}$. 
Let $L$ be a linear form not vanishing on the single point of $X-Z$.  Note
that $G$ also does not vanish on the single point of $X-Z$.  $LG$ is
certainly in the image of $\mu_{l,Z}$, but not in the image of $\mu_{l,X}$
because $LG$ does not vanish on all of $X$.

Case 2.
The image of $\mu_{l,X}$ has codimension one in $\ix_{l+1}$.  We need to show
that there is at least one subset ${\Bbb Z}$, with cardinality $d-1$ and
truncated Hilbert function, so that $\mu_{l,{\Bbb Z}}$ is surjective.  For this
case  we consider all subsets of
$X$ of cardinality
$d-1$.  Set
$Z_i = X-\{P_i\}$,
$i=1, \dots, d$.  Let $G_i$ be a form of degree $l$ in $I(Z_i)_l - \ix_l$. 
Note that $G_i$ is well defined up to elements of $\ix_l$.  One can see that
$F_1, \dots, F_t, G_1, \dots, G_d$ form a basis for $S_l$.  Indeed, since
$\ix_l$ has codimension $d$ in $S_l$ there are the right number of them to be
a basis, and any linear relation $a_1F_1 + \dots + a_tF_t + a_{t+1}G_1 +
\dots + a_{t+d}G_d =0$ would need to have $a_{t+i} = 0$, $i = 1, \dots, d$,
since $G_i(P_i)\ne 0$ but all the other $G$'s and $F$'s vanish at $P_i$. 
This would give a linear relation among the $F$'s which is impossible because
they form a basis for $\ix_l$.

We only need to find one $Z_i$ such that $\mu_{l,Z_i}$ is surjective, and
such that $Z_i$ has truncated Hilbert function.  We will first argue that in
this situation, if $\mu_{l,Z_i}$ is surjective then $Z_i$ {\em must} have
truncated Hilbert function.

Let $L$ be a general linear form and let 
\[
J = \frac{I(X)+(L)}{(L)} \hbox{\hskip 2cm} J_i = \frac{I(Z_i)+(L)}{(L)}
\]
be the corresponding ideals in $R = S/(L) \cong k[x,y]$.  Note that $L$ is
not a zero divisor on $S/I(X)$ or $S/I(Z_i)$.  By slight abuse of notation, we
will call the rings $A = R/J$ and $A_i = R/J_i$  the {\em Artinian reductions}
of $X$ and $Z_i$, respectively.  

\begin{clm}\label{same-gens}
$I(X)$ and $J$ (resp.\ $I(Z_i)$ and $J_i$) have the same number of minimal
generators, occurring in the same degrees.
\end{clm}
\begin{proof}
This is standard. See for instance \cite{birkhauser}, p.\ 28.
\end{proof}

\begin{clm}\label{case2-artinian}
$J_i$ is {\em equal} to $J$ in degrees greater than or equal to $l$ if and only
if $Z_i$ does {\em not} have truncated Hilbert function.
\end{clm}
\begin{proof}
Notice that $J \subset J_i$ for all $i$, and notice that $J_i = J$ in degrees
$\geq l+1$, so it is enough to prove that $\dim J_i =
\dim J$ in degree $l$ if and only if $Z_i$ does not have truncated Hilbert
function.

Consider the Hilbert functions of $A$ and of $A_i$:
\[
\begin{array}{ll}
h_A: & 1 \ \ a_1 \ \ a_2 \ \ \dots \ \ a_{l-1} \ \ a_l \ \ 0 \\
h_{A_i}: & 1 \ \ b_1 \ \ b_2 \ \ \dots \ \ b_{l-1} \ \ b_l \ \ 0 
\end{array}
\]
Since $J \subset J_i$ we have $a_j \geq b_j \geq 0$ for all $j$.  We also have
$\sum a_j = d$ and $\sum b_j = d-1$.  It follows that for {\em one} value of
$j$, say $j_0$, we have $a_{j_0} = b_{j_0} + 1$, and for all other $j$ we have
$a_j = b_j$.  Since $Z_i$ has truncated Hilbert function if and only if $j_0 =
l$, this completes the proof of the claim.
\end{proof}

It follows from Claims \ref{same-gens} and \ref{case2-artinian} that if $Z_i$
does not have truncated Hilbert function then it is impossible that $X$ has
a minimal generator in degree $l+1$ but $Z_i$ does not have a minimal generator
in degree $l+1$.  So if we prove the existence of a $Z_i$ with no minimal
generator in degree $l+1$ then the truncated Hilbert function will follow
automatically.

Suppose that $\mu_{l,Z_i}$ is never surjective.  Since
$\dim\so I(Z_i)_l =
\dim\so\ix_l +3$, $\dim I(Z_i)_{l+1} = \dim\ix_{l+1}+1$, and by assumption
\[
\dim\ix_{l+1} =  \dim\mu_{l,X}(\so\ix_l)+1,
\]
 we see that for every
$i$ the kernel of $\mu_{l,Z_i}$ must have dimension at least two larger
than the dimension of the kernel of $\mu_{l,X}$.  That is, there must be two
degree one relations of the form
\begin{align*}
L_{i, 1}F_1 + \dots + L_{i, t}F_t + L_{i, t+1}G_i &= 0\\
M_{i, 1}F_1 + \dots + M_{i, t}F_t + M_{i, t+1}G_i &= 0.
\end{align*}
These relations must be linearly independent of each other and no linear 
combination of the two of them can involve only $F$'s and not $G_i$.  From
this one can see that all $2d$ of these relations (as you vary $i$) are
linearly independent elements of the kernel of the multiplication map $\so
S_l \to S_{l+1}$ which remain independent modulo the kernel of $\mu_{l,X}$.

Using our assumption on the codimension of the image of $\mu_{l,X}$ in
$\ix_{l+1}$  we conclude that this image has codimension $d+1$ in
$S_{l+1}$.  Comparing $\mu_{l,X}$ with the multiplication map $\so S_l\to
S_{l+1}$ we see that $\dim \so S_l = \dim \so \ix_l +3d$.  However, from the
previous paragraph we know that the dimension of the kernel of $\so S_l \to
S_{l+1}$ is at least $2d$ larger than the dimension of the kernel of
$\mu_{l,X}$.  Counting dimensions we get that $\so S_l \to S_{l+1}$ is not
surjective.  But, it is a well known triviality that $\so S_l \to S_{l+1}$ is
surjective.  This contradiction finishes case 2.

Case 3.
The image of $\mu_{l,X}$ has codimension two in $\ix_{l+1}$.  This will 
involve quite a bit more work and will be done in section 5.

Case 4. 
The image of $\mu_{l,X}$ has codimension $c\ge 3$ in $\ix_{l+1}$.  Let 
$Z_i$, $F_i$, $G_i$, $J$ and $J_i$ be as in case 2.  The codimension of the
image of 
$\mu_{l,X}$ in $I(Z_i)_{l+1}$ is $c+1\ge 4$.  We want to show that there
is a
$Z_i$ with truncated Hilbert function, such that $I(Z_i)$ has $c-3$ minimal
generators in degree $l+1$.  By Claim \ref{case2-artinian}, if $Z_i$ does not
have truncated Hilbert function then $J_i = J$ in degrees $\geq l$.  Hence $J$
and $J_i$ have the same number of minimal generators in degree $l+1$, and by
Claim \ref{same-gens}, the same is true of $Z_i$ and $X$.  So just as in case
2, it is enough to prove the existence of a $Z_i$ with the right number of
minimal generators, and it will automatically have truncated Hilbert function.

The proof will be by induction on $d$.  Hence we can assume that the theorem is
true for all the $Z_i$, but suppose that it fails for $X$.  In this case we
conclude that for each $i=1, \dots, d$ we have at least one degree one relation
of the form $L_{i,1}F_1 + \dots + L_{i, t}F_t + L_{i, t+1}G_i = 0$.  If there
were always two or more such relations we could arrive at a contradiction as in
case 2, so assume for $i=1$ there is only one such relation.

As indicated above, we may assume that the theorem holds for $Z_1$.  Thus we
can find $P_j$, $j\in\{2, \dots, d\}$ such that $Z_{1, j} = Z_1 - \{P_j\}$
satisfies the conclusion of the theorem with respect to $Z_1$.  A basis for
$I(Z_1)_l$ consists of $F_1,
\dots, F_t, G_1$ and a basis for $I(Z_{1, j})_l$ consists of $F_1, \dots,
F_t, G_1, G_j$.  The relations $L_{i, 1}F_1 + \dots + L_{i, t}F_t + L_{i,
t+1}G_i = 0$ for $i=1, j$ say that the codimension of the
image of $\mu_{l,Z_1}$ in $I(Z_1)_{l+1}$ is exactly $c-1\ge 2$
(because we assumed only one such relation) and the codimension of the image of 
$\mu_{l,Z_{1, j}}$ in $I(Z_{1, j})_{l+1}$ is at least $c-2\ge
1$.  But the assumption that $Z_{1, j}\subset Z_1$ satisfies the theorem says
that the codimension of the image of 
$\mu_{l,Z_{1, j}}$ in $I(Z_{1, j})_{l+1}$ is $c-3$.  This
contradiction finishes case 4.
\end{proof}

\begin{rem} \label{wrong-approach}
The proof of Lemma \ref{lem:1} is surprisingly simple.  The idea is to start
with a subset $Y'$ of $X$ (beginning with any single point) and add one point of
$X$ at a time in such a way that at each step, the new subset
$Y$ has the predicted Hilbert function.  This is done by considering the linear
system of hypersurfaces of any degree $t$ containing $Y'$.  If the general
element of this linear system vanishes on all of $X$ then consider degree
$t+1$.  If not, there is some point $P$ of $X$ not in the base locus of this
linear system, and we take $Y = Y' \cup P$. 

One would naturally wonder if the same approach, building up to $X$ point by
point rather than taking point after point away from $X$, would similarly be an
easier approach to Theorem \ref{thm:4.6}.  In fact this seems to not work. 
Consider, for instance, a set of points with the following configuration:
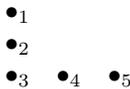
\[
\begin{array}{ccc}
\bullet_1 \\
\bullet_2 \\
\bullet_3 & \bullet_4 & \bullet_5
\end{array}
\]
(i.e.\ points 1, 2, and 3 are collinear and points 3, 4 and 5 are collinear). 
If we build up $X$ starting with point 3, it is of course possible to do so in
such a way that at each step the subset obtained has the right (truncated)
Hilbert function.  For example, the sequence 3, 1, 4, 2, 5 works.  However, it
is {\em impossible} to find a sequence beginning with 3 such that at each step
the subset has the right minimal free resolution according to Theorem
\ref{thm:4.6}.  Indeed, the only subset of $X$ consisting of four points and
having the right number of minimal generators is the set of points labeled
1, 2, 4 and 5.
\end{rem}

%%%%%%%%%%%%%%%%%%%%%%%%%%%%%%%%%%%%%%%%%%%%%%%%%%%%%%%%%%%%%%%%%%%%%%%%%%%%%

\section{The Final Case}

This section is devoted to proving case 3 of the proof of Theorem
\ref{thm:4.6}.  We are thus assuming that $X$ has two minimal generators in
degree $l+1$, and we are trying to show that there is a subset $Z_i$ of
cardinality $d-1$ having truncated Hilbert function and no minimal generator in
degree $l+1$.  With this assumption on $X$, the following fact is proved
exactly as in case 2 in the preceding section:  {\em  If a subset $Z_i$ of
$d-1$ points exists with no minimal generator in degree $l+1$ then it must have
truncated Hilbert function.}

 We begin by taking care of a special subcase.  For the following lemma we will
actually need to use the truncated Hilbert function to find the desired subset,
so we have to be a little careful.

\begin{lem}\label{gcd-case}
Assume that $X$ satisfies case 3, i.e.\ the image of $\mu_{l,X}$ has
codimension two in $I(X)_{l+1}$.  If the base locus of $I(X)_l$ is
one-dimensional then $X$ contains a subset, $Z$, of cardinality $d-1$ which
satisfies the rank condition asserted in Theorem \ref{thm:4.6}, namely $I(Z)$
has no minimal generators in degree $l+1$.
\end{lem}

\begin{proof}

By Lemma \ref{lem:1}, there is at least one subset $Z_i$
whose Hilbert function is the truncation of that of $X$.  We will find our
desired $Z$ from among these subsets, so from now on we will assume that this
is the Hilbert function of $Z$.  Then we have that the ideal of $X$ agrees
with that of $Z$ in degrees $\leq l-1$ and $Z$  and $X$ both impose
independent conditions on curves of degree $l$.  It follows that  
\[
\begin{array}{rcl}
\hbox{dim } I(X)_l + 1 & = & \hbox{dim } I(Z)_l, \\
\hbox{dim } I(X)_{l+1} + 1 & = & \hbox{dim } I(Z)_{l+1}
\end{array}
\]
We are assuming, furthermore, that $X$ has precisely two minimal generators in
degree $l+1$.  We need to show that $Z$ can be chosen with {\em no} minimal
generator in degree $l+1$.

The assumption about the dimension of the zero locus means that $I(X)_l$ has a
GCD, $F$.  Let $k$ be the degree of $F$.  By abuse of notation we will use
$F$ both for the curve in $\bp^2$ and for the polynomial.

\begin{clm}\label{clm1}
$k \leq 2$.
\end{clm}

\begin{proof}

We will use ideas from \cite{BGM} Proposition 2.3 (closely related to work of
Davis \cite{D}).  Let $X_1$ be the subset of $X$ lying on $F$ and let $X_2$ be
the subset not lying on $F$.  We have $I(X_1 ) = [I(X) +
(F)]^{\hbox{\footnotesize sat}}$ and $I(X_2) = [I(X) : F]$, which is already
saturated.  For $k \leq t
\leq l$ we have 
\begin{equation}\label{hfeqn}
\Delta h_{X_2} (t-k) = \Delta h_X (t) - k.
\end{equation}
Notice that $I(X)_l = F \cdot I(X_2)_{l-k}$.  From this we deduce two things. 
First, $X_2$ imposes independent conditions on forms of degree $l-k$ since $X$
does on forms of degree $l$. Second, $\hbox{rk } \mu_{l,X} = \hbox{rk }
\mu_{l-k,X_2}$.  

Let $J$ be the ideal generated by $I(X)_{\leq l}$.  We have just seen that
$\dim J_l = \dim I(X_2)_{l-k}$.  In degree $l+1$ we have the inequality $\dim
J_{l+1} \leq \dim I(X_2)_{l-k+1}$, where the failure to be an equality is
measured by the number of minimal generators of $I(X_2)$ in degree
$l-k+1$.  Let $h(S/J, t)$ be the corresponding Hilbert function and consider
$\Delta h(S/J,l+1)$.  From the above considerations, one can check that 
\[
\begin{array}{rcl}
\Delta h(S/J,l+1) & \geq & k + \Delta h_{X_2} (l-k+1) \\
& = & k.
\end{array}
\]
On the other hand, since $I(X)$ has two minimal generators in degree $l+1$, we
have $2 = \dim I(X)_{l+1} - \dim J_{l+1}$.  This gives
\[
\begin{array}{rcl}

k & \leq & \Delta h(S/J,l+1) \\
& = & l+2-\dim J_{l+1} + \dim J_l \\
& = & l+2 - \left [ \dim I(X)_{l+1} - 2 \right ] + \dim I(X)_l \\
& = & \Delta h_X (l+1) + 2 \\
& = & 2

\end{array}
\]
and this proves the claim.
\end{proof}

\begin{clm}
If $k=2$ then $X_1$ consists of exactly $2l+1$ points on $F$.  If $k=1$ then
$X_1$ consists of exactly $l+1$ points on $F$.
\end{clm}

\begin{proof}

Let us collect the following facts.

\begin{enumerate}

\item The initial degree of $I(X)$ is $\geq 2$.

\item $\Delta h_{X_2} (l-k+1) = 0$ since $X_2$ imposes independent conditions on
forms of degree $l-k$ ($k = 1,2$).

\item $\Delta h_X(l+1) = 0$.  

\item $\Delta h_X (l) \geq 2$.  This follows because we are assuming that $X$
has two minimal generators in degree $l+1$.  It can be seen, for example, by
applying  \cite{campanella} Theorem 2.1 (d), since in our situation
certainly $X$ is contained in a complete intersection of type $(\alpha, \beta)$
with $\alpha < \beta = l+1$.

\item $\deg X = \deg X_1 + \deg X_2$.

\item $\displaystyle \sum_t \Delta h_{X_2} (t) = \deg X_2$.

\item $\displaystyle \sum_t \Delta h_{X} (t) = \deg X$.

\end{enumerate}

If one now considers the Hilbert function of the Artinian reduction of $S/I(X)$
(i.e.\ the function given by $\Delta h_X (t)$) and applies the
equation (\ref{hfeqn}), in the case $k=2$ (resp.\ $k=1$) one gets from the above
facts that $\deg X_1 = 2l+1$ (resp.\ $\deg X_1 = l+1$) as claimed.
\end{proof}

We consider the cases $k=2$ and $k=1$ separately.  Suppose that $X_1$ consists
of $2l+1$ points on either a smooth conic or else a union of two lines.  In the
latter case, either one point lies at the intersection of the two lines or else
there are $l$ points on one line and $l+1$ points on the other.  (Otherwise $X$
fails to impose independent conditions on forms of degree $l$.)  In any of these
cases, the removal of a suitable point $P$ leaves $2l$ points which form the
complete intersection of $F$ and some curve $G$ of degree $l$.  Let $Z$ be the
subset of $X$ obtained by removing this point.  

We know that there exists an element of $I(Z)_l$ which is not in $I(X)_l$.  We
first claim that such an element must meet $F$ in finitely many points. 
Certainly the base locus of the linear system $|I(Z)_l|$ cannot contain all of
$F$ since then it contains the deleted point $P$, contradicting the fact that
$X$ imposes independent conditions on forms of degree $l$.  Hence the assertion
is clear if $F$ is irreducible.   Suppose that $F = L_1 L_2$ is reducible and
suppose (without loss of generality) that $L_1$ is in the base locus of
$|I(Z)_l|$.  Then any element of this linear system consists of the product of
$L_1$ with a homogeneous polynomial of degree $l-1$ containing the remaining
points of $X$.  By construction, the remaining points include $l$ points on
$L_2$, so in fact all of $F$ is in the base locus, a contradiction.

Hence without loss of generality we may assume that the subset of $Z$ lying on
$F$ is the complete intersection of $F$ and a form $G \in I(Z)_l$ (i.e.\ $G$
contains all of $Z$).  Now, if $(F_1,\dots,F_m)$ form a basis for
$I(X)_l$ and $(F_1,\dots,F_m,G)$ form a basis for $I(Z)_l$, then any linear
relation
\[
L_1 F_1 + \cdots + L_m F_m + LG = 0
\]
implies $L=0$ since no factor of $F$ is a factor of $G$, and $\deg F = 2$.  The
conclusion follows from this fact.

We now turn to the case $k=1$.  We have that $X$ is the union $X = X_1 \cup
X_2$, where $X_1$ consists of $l+1$ points on the line $F$.  Since
$|I(X)_{l+1}|$ has a zero-dimensional base locus, a general element, $G$, of
this linear system does not contain $F$ as a component.  Hence in particular
$X_1$ is the complete intersection of $F$ and $G$.  Also, in particular we have
that $G \in I(X_2)$.   We now note that $X = X_2 \cup X_1$ is a liaison
addition (cf.\ \cite{GM4}, \cite{schwartau})!  Hence its ideal is of the form
\[
I(X) = F \cdot I(X_2) + (G).
\]
Since $I(X)$ has two minimal generators of degree $l+1$, clearly $G$ must be
one of these and $I(X_2)$ must have exactly one minimal generator  in degree
$l$ (which is the maximum possible degree).

Let $Z$ be the subset of $X$ obtained by removing a point, $P$, of $X_1$.  Let
$Z_1$ be the subset of $X_1$ obtained by removing $P$. Exactly as above, $Z =
Z_1 \cup X_2$ is a liaison addition, since the linear system $|I(Z)_l|$ does
not have all of $F$ in its base locus.  Its ideal is of the form
\[
I(Z) = F \cdot I(X_2) + (G')
\]
where $\deg G' = l$, $G' \in I(X_2)$ and $(F,G')$ form a complete
intersection.  

We want to show that $P$ can be removed in such a way that $I(Z)$ has no
minimal generator in degree $l+1$.  

\begin{clm} \label{min-gen-claim}
The following are equivalent:
\begin{itemize}
\item[(1)] $I(Z)$ has a minimal generator in degree $l+1$.
\item[(2)] $I(X_2)$ has a minimal generator in degree $l$ other than $G'$.
\item[(3)] $G'$ is not a minimal generator for $I(X_2)$.

\end{itemize}
\end{clm}

\begin{proof}
The equivalence of (2) and (3) is clear since we have already observed that
$I(X_2)$ has exactly one minimal generator in degree $l$.  The fact that (1)
implies (2) follows from the equation $I(Z) = F \cdot I(X_2) + (G')$.  Now
assume (2), and let $H$ be the minimal generator in $I(X_2)$ other than $G'$. 
We want to show that $FH$ is a minimal generator for $I(Z)$.  Suppose not, and
let $F_1,\dots,F_k$ be a basis for $I(X_2)_{l-1}$.  Then we have
\[
FH = LG' + \sum_{i=1}^k L_i F F_i,
\]
or equivalently 
\[
F \left ( H-\sum_{i=1}^k L_i F_i \right ) = LG'.
\]
Since $F$ does not divide $G'$, we get that $H = G' + \sum L_i F_i$ up to
scalar multiples, contradicting the assumption that $H$ is a minimal generator
for $I(X_2)$ other than $G'$.
\end{proof}

As a result of Claim \ref{min-gen-claim} we have in particular that {\em $I(Z)$
has a minimal generator in degree $l+1$ if and only if $G'$ is {\em not} a
minimal generator of $I(X_2)$.}
We want to show that we can find the subset $Z$ with no minimal generator in
degree $l+1$.  We thus have to show that among the $l+1$ collinear points of
$G \cap F$, there is at least one point $P$ whose removal leads to a $G' \in
I(X_2)_l$ as above which is not in the image of $\mu_{l-1,X_2}$.  We will do
this by contradiction.  Let $P_1,\dots,P_{l+1}$ be the points of $G \cap F$. 
For each $i$, $1\leq i \leq l+1$, let $G_i \in I(X_2 \cup P_1 \cup \dots \cup
P_{i-1} \cup P_{i+1} \cup 
\dots \cup P_{l+1})_l \subset I(X_2)_l$ such that $G_i$ does not contain $F$ as
a factor, as was done with $G'$ above.  Note that $G_1,\dots,G_{l+1}$ are
linearly independent in $S_l$, as was done in Case 2 in the previous section.

We want to show that it is impossible for
$G_1,\dots,G_{l+1}$ to all be in the image of $\mu_{l-1,X_2}$.  Suppose
otherwise.  Note that $F$ is a non zero-divisor of $S/I(X_2)$, and let 
\[
J = \frac{I(X_2) + (F)}{(F)} \hskip 1.5cm \hbox{and} \hskip 1.5cm R = S/(F).
\]
We may view $R/J$ as the Artinian reduction of $S/I(X_2)$, and in particular
$J$ has a minimal generator in degree $l$ since we saw that $I(X_2)$ must have
a minimal generator in degree $l$.  On the other hand, if we let
$\bar G_i$ be the image of $G_i$ in $R$, the same argument as above gives that
the
$\bar G_i$ are linearly independent in $J_l \subset R_l$.  Note that $\dim J_l
= l+1$, so the $\bar G_i$ form a basis.  If the $G_i$ are all in the image of
$\mu_{l-1,X_2}$ then none of the $\bar G_i$ is a minimal generator of $J$, so
$J$ has no minimal generator in degree $l$, a contradiction.  This concludes
the proof of Lemma \ref{gcd-case}.
\end{proof}

Let $Z_i$, $F_i$, and $G_i$ be as in case 2.  As noted above, the truncated 
Hilbert function will follow immediately once we find the $Z_i$ with no minimal
generator in degree $l+1$.  Since we are assuming that the proposition fails
for $X$, we conclude that for each $i = 1, \dots, d$ we have at least one
degree one relation of the form
\begin{equation}\label{eq:5.1}
L_{i, 1}F_1 + \dots + L_{i, t}F_t + L_{i, t+1}G_i = 0
\end{equation}
where $L_{i,t+1} \neq 0$.  If there were always two or more such relations we
could arrive at a  contradiction as in case 2.  So for some $i$'s there must be
exactly one such relation.

As a result of Lemma \ref{gcd-case}, we may assume that the base locus of
$|I(X)_l|$ is zero-dimensional.  Consequently we may choose homogeneous
polynomials $H,K \in I(X)$ {\em each of degree $\leq l$} which form a regular
sequence, hence link $X$ to a zeroscheme $D$.  Note that we may choose $H$
and $K$ to be minimal generators of $I(X)$.  Without loss of generality, say
$\deg K \leq \deg H \leq l$.

\begin{clm}\label{clm:3}
The ideal $I(D)$ of $D$ contains a form of degree less than $\deg(K)$.
\end{clm}
\begin{proof}
We use the Cayley-Bacharach Theorem \cite{DGO} Thm.~3(b).  The first difference 
function for the Hilbert function of the complete intersection $H\cap K$ is:
\begin{equation*}
1, 2, 3, \dots, \deg K, \deg K, \dots, \deg K, \deg K-1, \deg K-2, \dots, 3, 
2, 1,
\end{equation*}
with $\deg K$ repeated $\deg H - \deg K + 1$ times.  The first difference 
function for the Hilbert function of $X$ does not reach $0$ until degree
$l+1$.  Since $\deg H\le l$, when we subtract these and read backwards to get
the first difference function for $D$ we see its maximum is less than $\deg
K$.
\end{proof}

Now we wish to use liaison theory to compare the resolutions of $\ix$ and 
$I(D)$.  We follow the presentation in \cite{CGO}.  Suppose the resolution for
$\ix$ is
\begin{equation*}
0\to\bigoplus_{i=1}^{e+1} S(-m_i)\, \stackrel{\varphi}{\longrightarrow} 
\bigoplus_{i=1}^{e+2} S(-d_i) \to \ix \to 0
\end{equation*}
with $d_1\geq d_2 \ge \dots \ge d_{e+2}$, $m_1 \geq m_2 \ge \dots \ge
m_{e+1}$.  By our assumptions on generators of $\ix$ in degree $l+1$ we have
that $d_1 = d_2 = l+1$, $d_3\le l$.  The map $\varphi$ is given by a matrix of
forms
\begin{equation*}
\mathcal A = 
	\begin{bmatrix}
	a_{11} & a_{12} & \dots & a_{1, e+1} & a_{1, e+2}\\
	a_{21} & a_{22} & \dots &&\\
	\vdots & \vdots &       & \vdots & \vdots \\
	a_{e+1,1} & a_{e+1,2} & \dots & a_{e+1, e+1} & a_{e+1, e+2}
	\end{bmatrix}
\end{equation*}
As in \cite{CGO} we use $\partial \mathcal A$ to denote the integer matrix 
whose $(i, j)$ entry is $\deg a_{ij} = u_{ij} = \max \{0,m_i - d_j \}$.
\begin{equation*}
\partial \mathcal A = 
	\begin{bmatrix}
	u_{11}  & \dots & u_{1, e+2}\\
	\vdots  & \vdots \\
	u_{e+1,1} & \dots & u_{e+1, e+2}
	\end{bmatrix}.
\end{equation*}
The fact that $\ix$ has generators in degree $l+1$ says that $X$ has 
defining equations of high degree as defined in \cite{CGO} after Prop.~1.2,
so that Cor.~2.5 of \cite{CGO} says $u_{1, 1} = 1$.  This  together with
inequalities found in Remark~2.2 of \cite{CGO} and the previously mentioned
facts about $d_1, d_2, d_3$ says that $\partial \mathcal A$ has the form
%%%%%%%%%%
%%%%%%%%%%
\begin{equation*}
\begin{bmatrix}
\vphantom{[} & & \\ \quad
	\begin{matrix}
	1  & \quad & 1\\
	1  & \quad & 1\\
	   & \vdots & \\
	1  &  & 1
	\end{matrix} 
\qquad & \vline &\quad \quad \quad  A \quad \quad  \qquad\\
\vphantom{[} &\vline & \\
\hline&\vline\\
	\begin{matrix}
	0  & \quad & 0\\
	   & \vdots & \\
	0  &  & 0
	\end{matrix} 
& \vline & B\\
\vphantom{[} & & 
\end{bmatrix}
\end{equation*}
%%%%%%%%%%%%%%%
%%%%%%%%%%%%%%%%
\noindent where all entries of $A$ are greater than one.  Applying the 
description of liaison of \cite{CGO} section 3 we get the corresponding
matrix for $I(D)$ by eliminating from $\partial\mathcal A$ the columns
corresponding to $H$ and $K$, taking the transpose, and reversing the order
of columns and rows to get something that looks like
%%%%%%%%%%
%%%%%%%%%%
\begin{equation*}
\begin{bmatrix}\quad 
\qquad & \vline &\\
 B'& \vline &A'\\
\vphantom{[} &\vline &  \\
\hline&\vline&\:\\ \quad
	\begin{matrix}
	0  & \quad & 0\\
	   & \hdots & \\
	0  &  & 0
	\end{matrix} 
&\vline
&	\begin{matrix}
	1  & \quad & 1\\
	   & \hdots & \\
	1  &  & 1
	\end{matrix} 
\vphantom{[} & & 
\end{bmatrix}
\end{equation*}
%%%%%%%%%%%%%%%
%%%%%%%%%%%%%%%%
\noindent where all the entries of $A'$ are greater than 1.  One way to look 
at this is to say that the two minimal generators of $\ix$ of maximal degree
$l+1$ induce two degree 1 relations among the forms of lowest degree in
$I(D)$.

We can also use $H$ and $K$ to link $Z_i$ and $D\cup \{P_i\}$ for each $i = 
1, 2, \dots, d$.  (Here and subsequently, we abuse notation somewhat and write
$D \cup \{ P_i \}$ for the scheme residual to $Z_i$, even when $P_i$ is a
common component of $X$ and $D$.)  Unfortunately, $H$ and $K$ may not be
minimal generators for $Z_i$, so we have to argue slightly differently to get
the desired matrix for $I(D \cup \{P_i\})$.  We have seen that if there is a
$Z_i$ with no minimal generator in degree $l+1$ then it has truncated Hilbert
function, and it is the subset that we are seeking.  Hence without loss of
generality we may assume that $Z_i$ has at least one minimal generator in
degree $l+1$, and, as above, it has a degree matrix of the form
\begin{equation*}
\begin{bmatrix}
\vphantom{[} & & \\ \quad
	\begin{matrix}
	&1  \\
	& \vdots \\
	&1
	\end{matrix} 
\qquad & \vline &\quad \quad \quad A_i \quad \quad \qquad\\
\vphantom{[} &\vline & \\
\hline&\vline\\
	\begin{matrix}
	& 0\\
	& \vdots & \\
	& 0
	\end{matrix} 
& \vline & B_i\\
\vphantom{[} & & 
\end{bmatrix}
\end{equation*}
\noindent where all entries of $A_i$ are greater than or equal to one. As a
result, the minimal free resolution for $I(Z_i)$ has the form
\begin{equation*}
0\to\bigoplus_{i=1}^{f+1} S(-n_i)\, \stackrel{\varphi}{\longrightarrow} 
\bigoplus_{i=1}^{f+2} S(-e_i) \to I(Z_i) \to 0
\end{equation*}
with $l+1=e_1\geq e_2 \ge \dots \ge e_{f+2}$, $l+2= n_1 \geq n_2 \ge \dots \ge
n_{f+1}$.
(This can also be deduced directly by considering the regularity of $I(Z_i )$.)
Now we link using $K$ and $H$ with $k = \deg K \leq h = \deg H < l+1$.  Let $C$
be the complete intersection of $H$ and $K$.  We get a commutative diagram
\[
\begin{array}{ccccccccccccc}
0 & \to &\displaystyle \bigoplus_{i=1}^{f+1} S(-n_i) &
\stackrel{\varphi}{\longrightarrow}  &\displaystyle \bigoplus_{i=1}^{f+2}
S(-e_i) & \to & I(Z_i) & \to & 0 \\
&& \uparrow && \uparrow && \uparrow \\
0 & \to & S(-h-k) & \to & S(-h) \oplus S(-k) & \to & I(C) & \to & 0
\end{array}
\]
and by the usual mapping cone trick (cf.\ \cite{birkhauser} Proposition 5.2.10)
we get a resolution for $I(D \cup \{P_i\})$ of the form
\[
0 \to \bigoplus_{i=1}^{f+2} S(e_i -h-k) \to \bigoplus_{i=1}^{f+1} S(n_i -h-k)
\oplus S(-h) \oplus S(-k) \to I(D \cup \{P_i\}) \to 0.
\]
This is not necessarily minimal.  It may be that zero, one or two terms can be
split off, depending on whether neither, one or both of $H$ and $K$ are minimal
generators of $I(Z_i)$, respectively.  However, by the assumption that $k \leq
h < l+1$, we see that in any case the smallest term $S(e_1-h-k)$ is {\em not}
split off.  It follows that the degree matrix for $I(D \cup \{P_i\})$ is of the
form
%%%%%%%%%%
%%%%%%%%%%
\begin{equation*}
\begin{bmatrix}\quad 
\qquad & \vline &\\
 B_i'& \vline &A_i'\\
\vphantom{[} &\vline &  \\
\hline&\vline&\:\\ \quad
0 \cdots 0 &\vline  &1 \cdots 1 \quad\\ 
\vphantom{[} & & 
\end{bmatrix}
\end{equation*}
%%%%%%%%%%%%%%%
%%%%%%%%%%%%%%%%
\noindent where all the entries of $A_i'$ are greater than or equal to one.  In
this  case we can say that the  minimal generator of $I(Z_i)$ of maximal degree
$l+1$ induces a degree 1 relation among the forms of lowest degree in
$I(D\cup\{P_i \})$.

Putting these facts together we get the following.  Let $v$ be the lowest 
degree of forms in $I(D)$.  We have two independent degree 1 relations among
$I(D)_v$.  For each $i = 1, 2, \dots, d$ there exists a linear combination
(depending on $i$) of these two relations that becomes a degree 1 relation on
$I(D\cup \{P_i\})_v$.

We will convert this to a question about sections of twisted bundles of 
differential forms on $\bp^2$.  The basic idea is fairly standard; see for
instance \cite{HS} section 9 or \cite{B} section 1.  Let
$\Omega^1_{\bp^2}$ represent the sheaf of differential one-forms on $\bp^2$. 
Consider the dual of the Euler sequence twisted by $v+1$ (cf.\ \cite{H2}
Theorem II.8.13):
\begin{equation}\label{euler-seq}
0 \rightarrow \Omega_{\bp^2}^1 (v+1) \rightarrow {\mathcal O}_{\bp^2}
(v)^{\oplus 3} \rightarrow {\mathcal O}_{\bp^2} (v+1) \rightarrow 0
\end{equation}
Taking cohomology, we obtain
\[
0 \rightarrow H^0(\bp^2,\Omega_{\bp^2}^1 (v+1)) \rightarrow H^0(\bp^2,{\mathcal
O}_{\bp^2} (v)^{\oplus 3}) \stackrel{\alpha}{\longrightarrow}
H^0(\bp^2,{\mathcal O}_{\bp^2}(v+1)) \rightarrow \cdots.
\]
The map $\alpha$ is given by $\alpha(f_1,f_2,f_3) = x_0 f_1 + x_1f_2 +
x_2f_3$.  A degree 1 relation among elements of $I(D)_v$ can be written as
$x_0g_1 + x_2g_2 + x_2g_3 = 0$ where $g_i \in I(D)_v$ and therefore represents
an element of the kernel of $\alpha$ vanishing on $D$.  By exactness it is the
image of an element of $H^0(\bp^2,\Omega^1_{\bp^2} (v+1))$ vanishing on $D$. 
The same argument applies to relations on $I(D\cup \{P_i\})_v$.  We deduce that
the two degree one relations among $I(D)_v$ correspond to two linearly
independent sections, call them $s_1$ and $s_2$, in $H^0(\bp^2,\Omega^1_{\bp^2}
(v+1))$ that vanish on $D$, and the degree one relation on $I(D \cup \{P_i\})_v$
corresponds to a linear combination of $s_1$ and $s_2$ that also vanishes at
$P_i$.  (When $P_i \in \hbox{Support}(D)$ this must be interpreted
appropriately in terms of ideals.)  The locus of points where sections of
vector bundles fail to be linearly independent can be analyzed using Chern
classes.

Set $Y = \{ \  P \in \bp^2 \ | \ s_1(P) \hbox{ and } s_2(P) \hbox{ are linearly
dependent } \}$.  Note that since $s_1$ and $s_2$ are linearly independent as
elements of $H^0(\bp^2,\Omega^1_{\bp^2} (v+1))$, $Y \neq \bp^2$.  From \cite{F}
Example 14.3.2 we see that the class of $Y$ with an appropriate scheme
structure is the first Chern class of $\Omega_{\bp^2}^1 (v+1)$.  Another way to
see this is to observe that $s_1 \wedge s_2$ is a nonzero section of the line
bundle $\bigwedge^2 (\Omega_{\bp^2}^1 (v+1))$.  From the Whitney sum formula
\cite{F} Theorem 3.2(e) applied to the Euler sequence (\ref{euler-seq}) we get
that the first Chern class of $\Omega_{\bp^2}^1 (v+1)$ is $2v-1$ times a line. 
In other words, $Y$ is a curve of degree $2v-1$.  Certainly $Y$ passes through
$D$ because $s_1$ and $s_2$ vanish at $D$, and $Y$ passes through $X = \{
P_1,\dots,P_d \}$ because for each $i$ some linear combination of $s_1$ and
$s_2$ vanishes at $P_i$.

We need to know more about the local nature of $Y$ near the points of the
complete intersection $H \cap K$, which we denote by $C$.  By abuse of notation
we also let $Y$ represent a polynomial generating the homogeneous ideal of
$Y$.  We consider the ideals $I(X)$, $I(D)$ and $I(C)$.  Let $Q$ be a point of
$C$.  The subscript $Q$ on an ideal or polynomial means we are considering the
corresponding ideal or function in the local ring of $\bp^2$ at $Q$.  

\begin{itemize}

\item[(a)] Suppose $Q \in X$, $Q \notin \hbox{Support}(D)$.  Then $Y_Q \in
I(X)_Q = I(C)_Q$.

\item[(b)] Suppose $Q \in \hbox{Support}(D)$, $Q \notin X$.  Choose
local coordinates $x,y$ on $\bp^2$ centered at $Q$ and trivialize
$\Omega_{\bp^2}^1 (v+1)$ locally near $Q$.  Then each $s_i$ is given locally by
a pair of functions $s_i(x,y) = (s_{i,1}(x,y),s_{i,2}(x,y))$, $i=1,2$.  Saying
that $s_i$ vanishes on $D$ is saying that $s_{i,j}(x,y) \in I(D)_Q$ for both
$j=1$ and $j=2$.  The local equation of $Y$ near $Q$ is the determinant
\begin{equation*}
\biggl| \begin{matrix}
s_{1, 1}(x,y) & s_{1, 2}(x, y)\\
s_{2, 1}(x,y) & s_{2, 2}(x, y)
\end{matrix} \biggr|
= s_{11}s_{22} - s_{12}s_{21}.
\end{equation*}
Thus $Y_Q \in I(D)_Q^2 \subset I(C)_Q$ since $I(D)_Q = I(C)_Q$.

\item[(c)] Suppose $Q \in X \cap \hbox{Support}(D)$.  We continue with the
notation of (b).  Some linear combination of $s_1$ and $s_2$ vanishes on the
residual scheme to $X-Q$ in $C$, which we called $D \cup \{ Q \}$ by abuse
of notation.  In the determinant replacing $s_1$ by this linear combination and
$s_2$ by some other independent linear combination we get that $Y_Q \in I(C)_Q
\cdot I(D)_Q \subset I(C)_Q$.

\end{itemize}

Putting the three local calculations (a), (b) and (c) together, we get that
globally $Y \in I(C)$, so we may write $Y = SH+TK$.  Again consider $Q$ a point
of $\hbox{Support}(D)$.  Whether we are in case (b) or (c), we have $Y_Q \in
I(C)_Q \cdot I(D)_Q$.  Because $I(C)_Q = \langle H_Q,K_Q \rangle$, we may write
$Y_Q = aH_Q + bK_Q$, where $a,b \in I(D)_Q$.  On the other hand, of course we
also have $Y_Q = S_Q H_Q + T_Q K_Q$.  This gives $(S_Q -a)H_Q = (b-T_Q)K_Q$. 
The local ring of $\bp^2$ at $Q$ is a unique factorization domain, and $H_Q$
and $K_Q$ have no common factors.  We conclude that $K_Q$ divides $S_Q -a$, so
that $S_Q -a \in I(C)_Q \subset I(D)_Q$, giving $S_Q \in I(D)_Q$.  Similarly,
$T_Q \in I(D)_Q$.  Since this is true locally for all $Q \in
\hbox{Support}(D)$, we conclude that globally $S,T \in I(D)$.  Furthermore,
since $Y \neq \bp^2$, $S$ and $T$ cannot both be zero.

 From Claim \ref{clm:3}, we see that $H$ and $K$ both have degree greater than
$v =$ smallest degree of a form in $I(D)$.  Since $Y$ has degree $2v-1$, both
$S$ and $T$ have degree less than $v-1$, meaning that they could not be in
$I(D)$.  This final contradiction completes the proof of Theorem \ref{thm:4.6}. 
\qed

\begin{rem} \label{MRC-in-P2}
It should be noted that Theorem \ref{thm:4.6} gives another proof of the
Minimal Resolution conjecture for general sets of points in $\bp^2$.  (Other
solutions can be found in \cite{GGR}, \cite{Ger-Mar}, \cite{GM4}.  See also
\cite{L2}.)

To see this first note that by \cite{Ger-Mar} Prop.\ 1.4, any set, $\Bbb X$, of
$\binom{d+2}{2}$ points in $\bp^2$ not lying on a curve of degree
$d$, always has resolution
$$
0 \rightarrow S( -(d+2))^{d+1} \longrightarrow S( -(d+1))^{d+2}
\longrightarrow I({\Bbb X}) \rightarrow 0.
$$
Thus, if $\Bbb Z$ is a general set of $t = \binom{d+1}{2} + r$ points in
$\bp^2$, $0<r<d+1$, then $\Bbb Z$ satisfies the minimal resolution
conjecture if and only if the multiplication map
$$
\mu_{d,{\Bbb Z}} : S_1 \otimes I({\Bbb Z})_d \longrightarrow I({\Bbb Z})_{d+1}
$$
has maximal rank (see e.g.\ \cite{Ger-Mar} section 2), i.e.
\begin{equation}\label{tony-*}
\hbox{rank } \mu _{d,{\Bbb Z}} = \min \{ \ \dim I({\Bbb Z})_{d+1}, 3\dim I({\Bbb
Z})_d \ \}.
\end{equation}

Now, (\ref{tony-*}) follows immediately from Theorem \ref{thm:4.6} once we note
that since \break $I({\Bbb X})_d = 0$ then both $\hbox{rank } \mu_{d, {\Bbb X}}
= 0$ and $\dim ( S_1 \otimes I({\Bbb X})_d) = 0$.
\end{rem}

\end{document}